\def\bbuildrel#1_#2^#3{\mathrel{\mathop{\kern 0pt#1}\limits_{#2}^{#3}}}

\def\NN{\mathbb N}

\def\QQ{\mathbb Q}

\def\Pot{\hbox{$\mathcal P$}}

\def\bs{\bigskip}
\def\ms{\medskip}
\def\ss{\smallskip}

\def\w{\thinspace\hbox{\hsize 14pt \rightarrowfill}\thinspace}

\def\0{\hbox{$\emptyset$}}

\def\M{\hbox{\rm MGR}}

\def\s{\hbox{$\sigma$}}

\def\sub{\subseteq}

\def\M{\hbox{$\mathcal M$}}

\def\F{\hbox{$\mathcal F$}}

\def\Ra{\Rightarrow}

\def\B{\mathscr{B}}

\def\M{\mathscr{M}}

\def\U{\mathscr{U}}

\documentclass[a4paper,12pt]{amsart}

\usepackage{amssymb,enumerate,mathrsfs}
\usepackage[utf8]{inputenc} 

\theoremstyle{plain}

\newcommand{\cantor}{2^{\NN}}
\newcommand{\baire}{\NN^{\NN}}
\newtheorem{theorem}{Theorem}[section]

\newtheorem{example}[theorem]{Example}

\newtheorem{proposition}[theorem]{Proposition}

\newtheorem{remark}[theorem]{Remark}
\newtheorem{cor}[theorem]{Corollary}
\newtheorem{claim}[theorem]{Claim}

\numberwithin{equation}{section}

\begin{document}

\title{Countably perfectly meager sets}

\author{R. Pol and P. Zakrzewski}
\address{Institute of Mathematics, University of Warsaw, ul. Banacha 2,
02-097 Warsaw, Poland}
\email{pol@mimuw.edu.pl, piotrzak@mimuw.edu.pl}

\subjclass[2010]{03E20, 03E15, 54E52}

\date{March 19, 2021}

\keywords{perfectly meager set, universally meager set, perfectly meager set in the transitive sense, $\lambda'$-set,  \s-ideal}

\begin{abstract}

We study a strengthening of the notion of a perfectly meager set.

We say that that a subset $A$ of   a perfect Polish 
space $X$ is
 countably perfectly meager in $X$, if for  every sequence of perfect subsets $\{P_n: n \in \NN\}$ of $X$, there exists an $F_\sigma$-set $F$ in $X$ such
that $A \sub F$ and $F\cap P_n$ is meager in $P_n$ for each $n$.

 We give various characterizations and examples of countably perfectly meager sets. We prove that   not  every universally meager set  is countably perfectly meager correcting an earlier result of Bartoszy\' nski.

\end{abstract}

\maketitle

\section{Introduction}\label{sec:1} 
Let us recall that a subset $A$ of   a perfect Polish
space $X$ is {\it universally meager} ($A\in
\mathbf{UM}$, see \cite{z-1}, \cite{z-2}, \cite{b-1}, \cite{b-2}), if  for every Borel
isomorphism $f$ between $X$ and 
any perfect Polish
space $Y$ the image of $A$ under $f$ is meager 
 in $Y$ (this class of sets was earlier
introduced and studied by Grzegorek \cite{g-1}, \cite{g-2}, \cite{g-3}
under the name of {\it absolutely
	of the first category} sets).

 Let us also recall  that $A$ is
{\it perfectly meager} ($A\in \mathbf{PM}$), if for all perfect
subsets $P$ of $X$, the set $A\cap P$ is meager in $P$. Clearly, $A\in \mathbf{PM}$ if and only if for every  perfect subset $P$ of $X$, there exists an $F_\sigma$-set $F$ in $X$ such
that $A \sub F$ and $F\cap P$ is meager in $P$ (cf. \cite[Theorem 6]{b-2}).

We shall say that $A$ is
{\it countably perfectly meager} ($A\in \mathbf{PM}_\sigma$), if for  every sequence of perfect subsets $\{P_n: n \in \NN\}$ of $X$, there exists an $F_\sigma$-set $F$ in $X$ such
that $A \sub F$ and $F\cap P_n$ is meager in $P_n$ for each $n$.

\ms

It follows directly from the definition that the class $\mathbf{PM}_\sigma$ is, exactly as the other two classes, a \s-ideal of subsets of the underlying space $X$ (shortly: a \s-ideal on $X$), i.e., it is hereditary, closed under taking countable unions and contains all singletons.

 One readily checks that $\mathbf{UM}\sub \mathbf{PM}$ and it is consistent that $\mathbf{UM}\subsetneq \mathbf{PM}$ but also that $\mathbf{UM}= \mathbf{PM}$ (see \cite{b-1}).
 
  Bartoszy\'nski \cite[Theorem 7, $(3)\Ra (2)\Ra (1)$]{b-2} proved that 
$\mathbf{PM}_\sigma\sub\mathbf{UM}$. Actually, it is in that paper where the property used by us to define the class $\mathbf{PM}_\sigma$ first appeared (without any specific name) 
and where it was claimed that this property  characterizes 
universally meager sets in the Cantor space $\cantor$.
 Unfortunately, 
there is a flaw in the part of the argument showing 
the inclusion  $\mathbf{UM}\sub\mathbf{PM}_\sigma$ (cf. \cite[Theorem 7, $(1)\Ra (3)$]{b-2}).

In fact, the following theorem  immediately implies  that it is consistent (in particular, true under {\bf CH}) that there exists a universally meager subset of $\cantor$ which is not countably perfectly meager in $\cantor$.

\begin{theorem}\label{main}
	Let $T$ be a subset of  $\cantor$ of cardinality $2^{\aleph_0}$. There exist a set $H\sub T\times \cantor$ intersecting each vertical section $\{t\}\times \cantor$, $t\in T$, in a singleton and a homeomorphic copy $E$ of $H$ in $\cantor$  which is not a $\mathbf{PM}_\sigma$-set  in $\cantor$. 
	In particular, $T$ is a continuous injective image of $E$.
\end{theorem}

Under the notation from Theorem \ref{main} it follows that if  $T$ is universally meager then $H$ is universally meager as well  ($\mathbf{UM}$ being closed with respect to
preimages under continuous injections (see \cite{z-1}) and so is its homeomorphic copy $E$.
 A refinement of this argument also shows that, at least consistently, in contrast to both $\mathbf{PM}$ and $\mathbf{UM}$ the class $\mathbf{PM}_\sigma$ is  not closed with respect to  homeomorphic images (see Theorem \ref{lambda}(3)). Consequently, unlike in the case of $\mathbf{PM}$ and $\mathbf{UM}$,
  the statement that 
 a subspace $A$ of a Polish space is countably perfectly meager makes sense 
 only if we specify a Polish space $X$ in which it is embedded. We shall therefore speak of {\sl countably perfectly meager sets in $X$} unless $X$ is clear from the context.

\ms

In Section \ref{sec:2} we present various characterizations and some examples of countably perfectly meager sets in $\cantor$. These include (for the definitions see Section \ref{sec:2-2}):
\begin{itemize}
	\item sets perfectly meager in the transitive sense, in particular:

	\begin{itemize}
		\item  $\gamma$-sets,
	\item strongly meager sets.
	\end{itemize}

	\item sets with the Hurewicz property and no perfect subsets,

	\item $\lambda'$-sets.

\end{itemize}

\ms

Section \ref{sec:3} is largely devoted to a proof of Theorem \ref{main}. 
We derive from it various examples of subsets of $\cantor$ which are universally meager but not countably perfectly meager in $\cantor$. 

In particular, if  there is a $\lambda$-set in $\cantor$ of cardinality of the continuum, then there is also one which is not countably perfectly meager in $\cantor$ (see Theorem \ref{lambda}(2)). 

Moreover, if  there exists a  $\lambda'$-set in $\cantor$ of cardinality of the continuum, then there is also one whose homeomorphic copy is not countably perfectly meager in $\cantor$  (see Theorem \ref{lambda}(3)). This strengthens the result of Sierpiński \cite{si} that $\lambda'$-property  is  not invariant under homeomorphisms.

We end Section \ref{sec:3} with a proof based on one of the characterizations of Section \ref{sec:2} that  the class $\mathbf{PM}_\sigma$ is closed under products in the sense that 
if $A$ and $B$ are $\mathbf{PM}_\sigma$-sets in perfect Polish spaces $X$ and $Y$, respectively, then $A\times B$ is a $\mathbf{PM}_\sigma$-set in $X\times Y$ (see Theorem \ref{product}).

\ms

In Section \ref{sec:4} we gather some additional comments.

 In Subsection \ref{sec:4.1} we present an example of a countably perfectly meager set in $\cantor$ which has neither the Hurewicz property nor $\lambda'$-property (see Example \ref{example 1}). We also give an example of a countably perfectly meager set in $\cantor$ which is not perfectly meager in the transitive sense (see  Example \ref{example 2}).

Subsection \ref{sec:4.2} contains remarks on some \s-ideals related to the classes $\mathbf{PM}$ and $\mathbf{PM}_\sigma$. 

\bs

\section{Characterizations and examples of countably perfectly meager sets}\label{sec:2} 
\ms 
\subsection{Characterizations of countably perfectly meager sets}\label{sec:2-1} \hfill\null

\smallskip

In this subsection $A$ is always a subset of a perfect (i.e., with no isolated points) Polish (i.e., a separable completely  metrizable) topological space $X$ and $\mathbf{PM}_\sigma$ is the family of all subsets of $X$ which are countably perfectly meager in $X$.

Let us recall that $A$ is an {\it $s_0$-set} if for every perfect (i.e., non-empty, closed and with no isolated points) set $P$ there is a  copy of the Cantor set $K\sub P$ with $K\cap A=\emptyset$. Clearly, every perfectly meager set has property $s_0$. 

\begin{theorem}\label{Cantor char}
	The following are equivalent:
	\begin{enumerate}
		\item $A\in \mathbf{PM}_\sigma$.
		\item For every sequence $\{K_n: n \in \NN\}$ of  copies of the Cantor set in $X$ there is an $F_{\sigma}$-set $F$ in $X$ such that $A\sub F$ and $K_m\not\subseteq F$ for each $m\in\NN$.
		\item For every sequence $\{K_n: n \in \NN\}$ of  copies of the Cantor set in $X$ there are closed sets $F_n$ in $X$ such that $A\sub \bigcup_n F_n$ and $K_m\not\subseteq F_n$ for each $m, n\in\NN$.
		\item For every sequence $\{K_n: n \in \NN\}$ of  pairwise disjoint copies of the Cantor set in $X$ there are closed sets $F_n$ in $X$ such that $A\sub \bigcup_n F_n$ and $K_m\not\subseteq F_n$ for each $m, n\in\NN$.
		\item $A$ is an $s_0$-set and for every sequence $\{K_n: n \in \NN\}$ of  pairwise disjoint and disjoint from $A$ copies of the Cantor set in $X$ there are closed sets $F_n$ in $X$ such that $A\sub \bigcup_n F_n$ and $K_m\not\subseteq F_n$ for each $m, n\in\NN$.
	\end{enumerate}
\end{theorem}

\begin{proof}
The implications 	$(1)\Ra (2)\Ra (3) \Ra (4)$ are  obvious.

	\smallskip
	
To prove that $(4) \Ra (5) $, it only suffices to show that if $(4)$ holds, then $A$ is an $s_0$-set. To see this, let us fix a  perfect set $P$ and let $K_0, K_1,\ldots$ be pairwise disjoint copies of the Cantor set in $P$ such that
 each non-empty relatively open set in $P$ contains some $K_n$. Now, by (4),  $A \sub \bigcup_n F_n$ for some closed sets $F_n$ in $X$ such that no $F_n$ covers any $K_m$. It follows that $F_n \cap P$ is nowhere dense in $P$ for each $n$ and so there exists a perfect set $K \sub P$ disjoint from $\bigcup_n F_n$. Then $K\cap A=\emptyset$ as well.

	\smallskip

To prove that $(5)\Ra (1)$, we  shall need the following 
simple observation.

\begin{claim}\label{disjoint Cantors}
	For  every sequence of perfect sets $\{P_n: n \in \NN\}$ in $X$, there is a sequence of pairwise disjoint Cantor sets $\{K_n: n \in \NN\}$ with $K_n\sub P_n$ for each $n$.
\end{claim}

 Indeed, first pick points $x_i\in P_i$ with $x_i\neq x_j$ for $i\neq j$, and then choose successively the Cantor sets $K_n$ in $P_n$ disjoint from $K_i$, for $i < n$, and $\{x_j: j>n\}$. 

\smallskip

Now let us assume (5)  and  let $\{P_n: n \in \NN\}$ be a sequence of perfect subsets of $X$. 

For each $n$ let $\{U^n_m: m \in \NN\}$ be a basis of non-empty relatively open subsets of $P_n$ and for each $m$ let us pick a Cantor set $K^n_m\sub U^n_m$. By the claim and the  fact that $A$ is an $s_0$-set we may assume that the sets $K^n_m$,  $n,m\in \NN$, are pairwise disjoint and disjoint from $A$.

By (5), there are closed sets $F_i$ in $X$ such that $A\sub \bigcup_i F_i$ and $K^n_m\not\subseteq F_i$ for each $m, n, i\in\NN$. Letting $F=\bigcup_i F_i$ we easily conclude that  $F\cap P_n$ is meager in $P_n$ for each $n$.

	\end{proof}

As a corollary let us formulate a characterization of countably perfectly meager sets stated and partially proved by Bartoszy\'nski in \cite{b-1}.

\begin{theorem}\label{ctble dense char}
	The following are equivalent:
	\begin{enumerate}
		\item $A\in \mathbf{PM}_\sigma$.
		\item For every sequence of countable dense-in-itself sets $\{A_n: n \in~\NN~\}$ there  are sets $B_n \sub A_n$ such that $\overline{A_n}=\overline{B_n}$ for each $n\in\NN$ and $\bigcup_n B_n$ is a $G_\delta$-set relative to $A \cup \bigcup_n A_n$.
		\end{enumerate}
\end{theorem}

\begin{proof}
	$(1)\Ra (2)$ was proved by Bartoszy\'nski \cite[Theorem 7, $(3)\Ra (2)$]{b-1}.
	
	\smallskip
	
	$(2)\Ra (1)$. To prove that $A\in \mathbf{PM}_\sigma$, we shall check that $A$ satisfies condition $(5)$ of Theorem \ref{Cantor char}. 
	
	First note that $A \in \mathbf{PM}$ (for a proof see \cite[Theorem 6, $(2)\Ra (1)$]{b-1}; in fact, Bartoszy\'nski \cite[Theorem 7, $(2)\Ra (1)$]{b-1} proved that $A\in \mathbf{UM}$) which implies that $A$ is an $s_0$-set. 
	
	Let $\{K_n: n \in \NN\}$ be a sequence of  (pairwise disjoint but this assumption is superfluous)  copies of the Cantor set disjoint from $A$. We shall show that there is an $F_\sigma$-set $F$ such that $A\sub F$ and $K_n \setminus F\neq\emptyset$ for each $n$. To this end, for each $n$ let us choose a countable dense subset $A_n$ of $K_n$.   
	By $(2)$,  there are sets $B_n \sub A_n$ such that $\overline{B_n}=K_n$ for each $n\in\NN$ and a $G_\delta$-subset $G$ of $X$ such that 
	$(A \cup \bigcup_n A_n)\cap G=\bigcup_n B_n$. Let $F=X\setminus G$. Then $A\sub F$ and for each $n$ we have $B_n\sub K_n\setminus F$, so $F$ is as required.
	
	\end{proof}

\begin{remark}
	The characterization of $\mathbf{PM}_\sigma$-sets from Theorem \ref{ctble dense char} should be compared with the following characterization of $\mathbf{PM}$-sets by  Bennett, Hosobuchi, and Lutzer \cite{b-h-l}. (for a short proof see \cite[Theorem 6]{b-1}):
	
	The following are equivalent:
	\begin{enumerate}
		\item $A\in \mathbf{PM}$.
		\item For every countable dense-in-itself set $A_0$ there  exists a set $B_0 \sub A$ such that $\overline{A_0}=\overline{B_0}$ and $B_0$ is a $G_\delta$-set relative to $A \cup A_0$.
	\end{enumerate}
	
	\end{remark}

Our last characterization of $\mathbf{PM}_\sigma$-sets is closely related to a characterization of $\mathbf{UM}$ sets (cf. Remark \ref{UM_Baire_char}).

\begin{theorem}\label{Baire sp char}
	The following are equivalent:
	\begin{enumerate}
		\item $A\in \mathbf{PM}_\sigma$.
		\item For every continuous bijection $f:\baire\w X$ 
		 there are closed sets $F_n$ in $X$ such that $A\sub \bigcup_n F_n$ and $f^{-1}(F_n)$ is nowhere dense in $\baire$ for each $n\in\NN$.
	\end{enumerate}
\end{theorem}

\begin{proof}
$(1)\Ra (2)$. Since open sets in $\baire$ are mapped by $f$ onto Borel sets, one can choose a sequence $\{B_n: n \in \NN\}$  of Borel subsets of $X$ such that $\{f^{-1}(B_n): n \in \NN\}$ is a basis of the topology of $\baire$. For each $n$ let us pick a Cantor set $K_n\sub B_n$. Since $A\in \mathbf{PM}_\sigma$,  there are closed sets $F_n$ in $X$ such that $A\sub \bigcup_n F_n$ and $K_m\not\subseteq F_n$ for each $m, n\in\NN$. It readily follows that $f^{-1}(F_n)$ is nowhere dense in $\baire$ for each $n\in\NN$.

$(2)\Ra (1)$. To prove that $A\in \mathbf{PM}_\sigma$, we shall check condition $(4)$ of Theorem \ref{Cantor char}.

 Let  $\{K_n: n \in \NN\}$ be a sequence of  pairwise disjoint Cantor sets in $X$.
 
  A key observation  is the 
following fact.

\begin{claim}\label{extending top}
There is a continuous bijection $f:\baire\w X$ 
 such that $f^{-1}(K_n)$ is open in $\baire$ for each $n\in\NN$.

\end{claim}
To prove the claim, let $\tau$ be the (perfect Polish) topology of $X$ and 
 let us first extend $\tau$ to the topology $\tau'$ whose basic open sets are elements of $\tau$ and relatively open subsets of $K_n$'s. More precisely, let $\tau_n$ be the topology generated by $\tau\cup\{K_n\}$ and then let $\tau'$ be the topology generated by $\bigcup_n\tau_n$. 
  This topology is Polish (cf. \cite[13.A]{ke}) and it is easy to check that it is also perfect. It follows (cf. \cite[7.15]{ke}) that there is a bijection $f:\baire\w X$ 
   which is continuous in the sense of $\tau'$ hence also in the sense of $\tau$.

\bs

Having proved the claim we may now apply the assumption about $A$ to find closed sets $F_n$ in $X$ such that $A\sub \bigcup_n F_n$ and $f^{-1}(F_n)$ is nowhere dense in $\baire$ for each $n\in\NN$. But $f^{-1}(K_m)$ being open in $\baire$ we conclude that $K_m\not\subseteq F_n$ for each $m, n\in\NN$. This shows that $A$ satisfies condition $(4)$ of Theorem \ref{Cantor char} completing the proof.

\end{proof}

\begin{remark}\label{UM_Baire_char}
The characterization of $\mathbf{PM}_\sigma$-sets from Theorem \ref{Baire sp char} should be compared with the following characterization of $\mathbf{UM}$-sets (cf. \cite[Theorem 2.4]{z-2}):

The following are equivalent:
\begin{enumerate}
	\item $A\in \mathbf{UM}$.
	\item For every continuous bijection $f:\baire\w X$ 
	 there are sets $F_n$ in $X$ such that $A\sub \bigcup_n F_n$ and $f^{-1}(F_n)$ is closed and nowhere dense in $\baire$ for each $n\in\NN$.
\end{enumerate}

In particular, in view of Theorem \ref{Baire sp char}, this gives another proof  of the inclusion $\mathbf{PM}_\sigma\sub \mathbf{UM}$.

\end{remark}

\subsection{Examples of countably perfectly meager sets}\label{sec:2-2} 
Let us recall that a subset $A$ of $\cantor$ is {\sl perfectly meager in the transitive sense} ($A\in \mathbf{AFC'}$, cf. \cite{n-s-w}, \cite{n-w-1} and \cite{w}) if  for every perfect subset  $P$ of $\cantor$, there exists an $F_\sigma$-set $F$ in $X$ such
that $A \sub F$ and $F\cap (P+t)$ is meager in $P+t$ for each $t\in \cantor$ or, equivalently (cf. \cite[Lemma 6]{n-w-1}), if for every sequence $\{K_n: n \in \NN\}$ of  copies of the Cantor set in $\cantor$ there are closed sets $F_n$ in $\cantor$ such that $A\sub \bigcup_n F_n$ and $K_m+t\not\subseteq F_n$ for each $m, n\in\NN$ and $t\in \cantor$. Combining this with Theorem \ref{Cantor char} we get the following result which somewhat strengthens the fact that $\mathbf{AFC'}\sub \mathbf{UM}$ established  by Nowik and Weiss \cite[Theorem 2]{n-w-1} by a similar argument.

\begin{theorem}
	Every subset of  $\cantor$ which is perfectly meager in the transitive sense is countably perfectly meager in $\cantor$. 
\end{theorem}

As a corollary we obtain a list of some classical classes of sets which being perfectly meager in the transitive sense are countably perfectly meager as well.

\begin{cor}\label{AFC' sets}
The following collections of subsets  of  $\cantor$  are countably perfectly meager in $\cantor$:
\begin{enumerate}
	\item meager-additive sets,
	\item  $\gamma$-sets,
	\item strongly meager sets,
	\item Sierpi\'nski sets.

\end{enumerate}	
\end{cor}

\begin{proof}

(1). See \cite[Proposition 6.6]{zind} 	
	
(2). See \cite{n-s-w}. This also follows from (1), since by \cite[Proposition 3.7]{n-w-2}, every $\gamma$-set is meager-additive. 

\smallskip

(3). See \cite[Theorem 9]{n-s-w}.

\smallskip

(4). This follows from (3), since Pawlikowski \cite{paw} proved that every Sierpi\'nski set is strongly meager.
\end{proof}	

The following result gives more examples of universally meager sets which are countably perfectly meager as well.

\smallskip

Let us recall that given a perfect Polish space $X$ a set $A\sub X$

\begin{itemize}
	\item  has {\it the  Hurewicz property}, if every continuous image of $A$ in $\baire$ is bounded in the ordering  $\leq^*$  of eventual domination,
\item is a  $\lambda'$-set in $X$ if every countable set $D\sub X$ is relatively $G_{\delta}$ in $A \cup D$.
\end{itemize}

The cardinal number $\mathfrak{b}$ is the minimal cardinality of a  subset of $\baire$ which is unbounded in the ordering  $\leq^*$.

\begin{proposition}\label{examples}
	The following collections of sets   are countably perfectly meager in the respective spaces:
	\begin{enumerate}
		\item subsets of a perfect Polish space $X$ with the Hurewicz property and no perfect subsets; in particular, subsets of $X$ of cardinality less than $\mathfrak{b}$,

		\item $\lambda'$-subsets in a perfect Polish space $X$, in particular:
	
		\begin{enumerate}
			
		\item sets in 
			$\baire$ of the form $\{f_\alpha:\alpha< \mathfrak{b}\}$ 
			 where
		
		\begin{itemize}
			\item $\alpha<\beta<\mathfrak{b}$ implies $f_\alpha<^*f_\beta$,
			\item for every $f\in\baire$ there is $\alpha<\mathfrak{b}$ with $f_\alpha\nleq^*f$.
		\end{itemize}

				\item Hausdorff $(\omega_1,\omega_1^{*})$-gaps in $\Pot(\NN)$.
			\end{enumerate}
	\end{enumerate}
	
\end{proposition}

\begin{proof}
(1). This was actually shown in Proposition 2.3 of Zakrzewski \cite{z-1}.	

\smallskip

(2). Let $A$ be a $\lambda'$-set in $X$. To prove that $A\in \mathbf{PM}_\sigma$, we shall check condition $(5)$ of Theorem \ref{Cantor char}. Clearly, $A$ is an $s_0$-set. For a sequence $\{K_n: n \in \NN\}$ of  pairwise disjoint and disjoint from $A$ copies of the Cantor set in $X$ for each $n$ let us pick a point $d_n\in K_n$ and let $D=\{d_n: n\in\NN\}$. Then, $A$ being a $\lambda'$-set, there is an $F_\sigma$ set $F$ in $X$ such that $A\sub F$ and $F\cap D=\emptyset$, so $d_n$ witnesses that  $K_n\not\subseteq F$ for any $n\in\NN$.

\smallskip

Sets described in (a) and (b) are classical examples of $\lambda'$ sets due to Rothberger and Hausdorff (see \cite{mi-1}).

\end{proof}	

\begin{remark}\label{Hurewicz vs lambda'}
	An easier way of proving that every Sierpi\'nski set in $\cantor$  in is
	$\mathbf{PM}_\sigma$ (cf. Corollary \ref{AFC' sets}) is to combine  Proposition \ref{examples}(1) with Theorem 7 of Fremlin and Miller \cite{f-m} which states that  every Sierpi\'nski set has the Hurewicz property. 
	
	Likewise, another way of proving that every $\gamma$-set in $\cantor$ is in 
	$\mathbf{PM}_\sigma$ (cf. Corollary \ref{AFC' sets}) is to  combine  Proposition \ref{examples}(1) with  Theorem 2 of Galvin and Miller \cite{g-m}	which states that every $\gamma$-set has the Hurewicz property. 
	
	On the other hand, the set described in Proposition \ref{examples}(2)(a), is $\lambda'$ in $\baire$ but does not have the Hurewicz property as an unbounded subset of $\baire$. Likewise, not every subset of $\cantor$ with the Hurewicz property and no perfect subsets (cf. Proposition \ref{examples}(1)) is a $\lambda'$-set in $\cantor$, see Example \ref{example 1} in the comment section below.

\end{remark}

\section{$\mathbf{PM}_\sigma$ versus $\mathbf{UM}$  }\label{sec:3}

Theorem \ref{main}, which we are now going to prove,  reveals an essential difference between the classes $\mathbf{UM}$ and $\mathbf{PM}_\sigma$.

\begin{proof}[Proof of Theorem \ref{main}]
	
	Let $C_0, C_1,\ldots$ be pairwise disjoint meager Cantor sets in $\cantor$ such that:

	\begin{enumerate}
		\item[(1)]  each non-empty open set in $\cantor$ contains some $C_n$. 
	\end{enumerate}
	
	\smallskip
	
	Let $P = \cantor\setminus\bigcup_n C_n$.
	
	\smallskip
	
	We shall justify the theorem in three steps (A), (B) and (C). 
	
	\medskip
	
	{\bf (A)} We claim that
	there exists a set $H\sub T\times P$ intersecting each vertical section $\{t\}\times P$, $t\in T$, in a singleton, such that each $F_\sigma$-set in $\cantor\times \cantor$ containing $H$ contains also $\{t\}\times V$ for some $t\in T$ and a non-empty open set $V$   in $\cantor$.

	Indeed, let $\{F_t: \ t\in T\}$ be a parametrization on $T$ of all  $F_\sigma$-sets in $\cantor\times \cantor$. For each $t\in T$, we pick $(t,\varphi(t))\in (\{t\}\times P)\setminus F_t$, whenever this is possible, and we let $\varphi(t)$ be an arbitrary fixed element of $P$, otherwise.
	
	Let us check that the graph $H=\{(t,\varphi(t)): t\in T \}$ has the required property.
	
	Let $F$ be an $F_\sigma$-set in $\cantor\times \cantor$ containing $H$, and let $t\in T$ be such that $F=F_t$. Then $(t,\varphi(t))\in F_t$, hence $F_t$ contains $\{t \}\times P$. Consequently, $P$ being a dense $G_\delta$-set in $\cantor$, the Baire category theorem provides a non-empty open set $V$ in $\cantor$ with $\{t \}\times V\sub F_t$, completing the proof of the claim.
	
	\medskip
	
	{\bf (B)} For any $s\in 2^{<\NN} $ let $N_s=\{x\in \cantor: s\sub x \}$
	be the standard 
	basic open set in $\cantor$ determined by $s$.

	Let $\sim$ be the equivalence relation on $\cantor\times \cantor$, whose equivalence classes are 
	given by:
	$$
	[(x,y)]_{\sim}=\left\{ \begin{array}{ll} N_{x|n}\times \{y\}, & \hbox{if}\ y\in C_n, \\ \{(x,y) \}, & \hbox{if}\  y\in P \end{array}\right.$$

	  Let $\pi(x,y)=[(x,y)]_{\sim}$ be the quotient map onto the quotient space $K=({\cantor\times\cantor})/\sim$ (whose topology consists of sets $U\sub K$ such that $\pi^{-1}(U)$ is open in $\cantor\times\cantor$).
	
	\medskip
	
	{\bf Claim.}  The space $K$ is  homeomorphic to $\cantor$.
	
	\medskip
	The equivalence classes of $\sim$ form an upper-semicontinuous decomposition of $\cantor\times \cantor$ (i.e., the saturation of every closed set in $\cantor\times\cantor$ is closed). 
	It follows that the decomposition space $K={\cantor\times\cantor}/\sim$ 
	is metrizable (cf. \cite[Theorem 4.2.13]{eng-2}). Moreover, $K$ is compact, zero-dimensional and has no isolated points, and hence the claim follows.
	However, for reader's convenience, we shall provide a direct argument to that effect, avoiding the metrization theorem.

	\medskip 
	
	Let $\B$ consist of clopen subsets of $\cantor\times \cantor$ of the form $N_s\times N_t$, where $s,t\in 2^{<\NN}$ and $N_t\cap C_k =\emptyset$ for each $k<{\rm length}(s)$. We shall show that $\{\pi(B):B\in\B\}$
	is a countable basis for $K$ consisting of clopen sets.
	
	\smallskip
	
	First, let us note that each set from $\B$ is saturated, i.e., is the union of equivalence classes. Indeed, if  
	$(x, y) \in N_s \times N_t$, $N_t\cap \bigcup_{k<n}C_k=\emptyset$ and $n={\rm length}(s)$, then 
	either $y\in P$ and then $[(x,y)]_{\sim}=\{(x,y)\}$ or $y\in C_m$ for some $m\geq n$ which implies that  $[(x,y)]_{\sim}=
	N_{x|m}\times \{y\}\sub N_s \times N_t$.
	
	It follows that for each $B\in \B$, $\pi(B)$ and $\pi((\cantor\times\cantor)\setminus B)$  are disjoint open sets in $K$, hence $\pi(B)$ is clopen in $K$.
	
		\smallskip
		
	Next, let us fix an open set $W$ in $K$ and let $c=\pi(x,y)\in W$. Then, since $\pi^{-1}(W)$ is open in $\cantor\times \cantor$ and 
	$(x,y)\in \pi^{-1}(W)$, we have	$N_{x|n}\times N_{y|m}\sub \pi^{-1}(W)$ for some $n$ and $m$. Moreover, if $y\in \bigcup_k C_k$, then we additionally assume that $n$ is the unique $k$ for which $y\in C_k$ (let us note that in this case $(x,y)\in N_{x|n}\times \{y\}\sub \pi^{-1}(W)$).   In any case, $y\notin\bigcup_{k<n} C_k$ and $\bigcup_{k<n} C_k$ being closed, there is large enough $m'\geq m$ for which $N_{y|m'}\cap \bigcup_{k<n} C_k=\emptyset$. Then
	$B=N_{x|n}\times N_{y|m'}\in \B$ and $\pi(B)$ is a neighbourhood of $c$ contained in $W$.
	
	\smallskip
	
	We have checked that $\{\pi(B):B\in\B\}$
	is a countable basis for $K$ consisting of clopen sets.
	Clearly, no $\pi(B)$ is a singleton, hence $K$ has no isolated points. 
	
	Finally,  the equivalence classes of $\sim$ are closed in
	$\cantor\times\cantor$, hence the singletons of $K$ are closed.

	It follows that the space $K$, being $T_1$ and having a basis consisting of clopen sets, is also Hausdorff and it is compact as a continuous image of $\cantor\times\cantor$. 
 Consequently,  being a compact, Hausdorff, second countable, zero-dimensional topological space without isolated points, $K$ is
	homeomorphic to $\cantor$, which completes the proof of the claim. 
	
	\medskip 
	
Let us also note that the sets	
  
	\begin{enumerate}
		\item[(2)] $P_s=\pi(N_s\times C_n)$, 
	\end{enumerate}
	where $s\in 2^{<\NN}$ and $n={\rm length}(s)$, are  perfect subsets of $K$. 
	
	\medskip
	
	{\bf (C)} 
	Finally, 
	let $E=\pi(H)$ (cf. (A)). Clearly, $E$ is a homeomorphic copy of $H$ in $K$ and $T$ is the injective image of $E$ under the continuous  function 
	$\hbox{proj}_1\circ \pi^{-1}|E$, where $\hbox{proj}_1$ is the projection of $\cantor\times\cantor$ onto the first axis.
	
	\smallskip
	
	We shall show that 
	
	\begin{enumerate}
		\item[(3)] $E$ is not a $\mathbf{PM}_\sigma$-set in $K$.
	\end{enumerate}

	To that end, let us consider an $F_\sigma$-set $F^*$ in $K$  such that $E\sub F^*$. Then 
	\begin{enumerate}
		\item[(4)] $F=\pi^{-1}(F^*)$ 
	\end{enumerate}
	is an $F_{\sigma}$-set in $\cantor\times\cantor$ containing $H$, so there are $t\in T$ and a non-empty open set $V$  in $\cantor$ such   $\{t\}\times V\sub F$, cf. (A).
	
	Let us fix $C_n\sub V$ (cf. (1)) and let $s=t|n$ be the unique sequence in $2^n$ such that 
	$t\in N_s$. We have $\{t\}\times C_n\sub F$ and let us notice that  $\pi(\{t\}\times C_n)=P_s$ (cf. (2)).
	
	Consequently, $P_s\sub F^*$, cf. (4).
	
	It follows that any $F_\sigma$-set in $K$ containing $E$ also contains some $P_s$, which confirms (3), completing the proof of the theorem.

\end{proof}	

As a corollary we have the following result which shows that, at least consistency--wise, the classes $\mathbf{UM}$ and  $\mathbf{PM}_\sigma$ are different (part (1)). 

Its part (2) strengthens the result of Nowik, Scheepers and Weiss \cite{n-s-w} that assuming the Continuum Hypothesis there is a $\lambda$-set in $\cantor$ which is not  perfectly meager in the transitive sense ($A$ a  {\it $\lambda$-set} if every countable set $D\sub A$ is relatively $G_{\delta}$ in $A$).

Part (3) strengthens the result of Sierpi\'nski \cite{si} that (assuming the continuum hypothesis) $\lambda'$-property  is  not invariant under homeomorphisms.

\begin{theorem}\label{lambda}\hfill\null
	\begin{enumerate}
		\item If  there exists a universally meager set in $\cantor$ of cardinality of the continuum, then there is also one which is not countably perfectly meager.
		
		\item If  there exists a  $\lambda$-set in $\cantor$ of cardinality of the continuum, then there is also one which is not countably perfectly meager.
		
		\item If  there exists a  $\lambda'$-set in $\cantor$ of cardinality of the continuum, then there is also one
			whose homeomorphic copy in $\cantor$  is not countably perfectly meager. In particular, the class $\mathbf{PM}_\sigma$ is  not closed with respect to  homeomorphic images.	
		
	\end{enumerate}
	In particular, 
	assuming the continuum hypothesis there is a $\lambda'$-set in $\cantor$ 
	whose homeomorphic copy in $\cantor$ is not countably perfectly meager.	
\end{theorem}

\begin{proof} We keep the notation from the proof of Theorem \ref{main}.	For a set $T\sub\cantor$ of cardinality $2^{\aleph_0}$ let 
	$H\sub T\times P$ and 
	$E\sub\cantor$ satisfy the assertions of Theorem \ref{main}.
	
	\smallskip
	
	(1) and (2).
	It suffices to notice that if $T$ is either universally meager or a $\lambda$-set, then so is $E$, respectively (cf. \cite{z-1}). But $E\not\in\mathbf{PM}_\sigma$.
	
	\smallskip
	
	(3). It can be readily checked that if $T$ is a $\lambda'$-set in $\cantor$, then $H$ is $\lambda'$-set in $\cantor\times P$. Since $P$ can be homeomorphically embedded as a dense $G_\delta$-set in $\cantor$, we may identify $H$ with a $\lambda'$-set in $\cantor$. But then $E$ is a homeomorphic copy of $H$ in $\cantor$ which is not $\mathbf{PM}_\sigma$ in $\cantor$.

\end{proof}

We close this section by showing that like $\mathbf{UM}$ (see \cite{z-1}) (but, at least consistency--wise, unlike $\mathbf{PM}$, see \cite{re-1}), the class $\mathbf{PM}_\sigma$ is closed under products.

\begin{theorem}\label{product}
	The product of two countably perfectly meager 
	sets is countably perfectly meager  in the sense that 
	if $A$ and $B$ are $\mathbf{PM}_\sigma$-sets in perfect Polish spaces $X$ and $Y$, respectively, then $A\times B$ is a $\mathbf{PM}_\sigma$-set in $X\times Y$.
\end{theorem}

\begin{proof}
	Let $A$ and $B$ be $\mathbf{PM}_\sigma$-sets in perfect Polish spaces $X$ and $Y$, respectively.
	
	To prove that $A\times B$ is a $\mathbf{PM}_\sigma$-set
	in $X\times Y$, we shall check condition (3) of Theorem \ref{Cantor char}.
	
	Let   $\{K_n: n \in \NN\}$ be a sequence of Cantor sets in $Z=X\times Y$ and for each $n$ let $L_n$ and $M_n$ be the images of $K_n$ under projections of $Z$ onto $X$ and $Y$, respectively. 
	
	Let $S=\{n\in\NN: L_n \hbox{ is uncountable} \}$ and for each $n\in S$ let us pick a Cantor set $L_n'\sub L_n$. Likewise, let  $T=\{n\in\NN: M_n \hbox{ is uncountable} \}$ and for each $n\in T$ let us pick a Cantor set $M_n'\sub M_n$.

	The sets $A$ and $B$ being countably perfectly meager, there are sequences $\{F^A_i:i\in \NN\}$ and $\{F^B_j:j\in \NN\}$ of closed sets  in $X$ and $Y$, respectively, such that $A\sub \bigcup_i F^A_i$, 
	$B\sub \bigcup_j F^B_j$, $L_m'\not\subseteq F^A_i$   and
	$M_n'\not\subseteq F^B_j$ whenever $m\in S,\ n\in T$ and $n,m,j,i\in\NN$.
	
	Let $F_{i,j}= F^A_i \times F^B_j$ for $i,j \in\NN$.
	
	Clearly, $A\times B \sub \bigcup_{i,j} F_{i,j}$ and we claim that for any $n, i, j$ we have $K_n\not\subseteq F_{i,j}$.
	
	Indeed, let us notice that $S\cup T =\NN$ since for each $n$ we have $K_n\sub L_n\times M_n$. It follows that  either $L_n'\not\subseteq F^A_i$ (if $n\in S$) or  $M_n'\not\subseteq F^B_j$ (if $n\in T$), so in either case $K_n\not\subseteq F_{i,j}$. This completes the proof.     
\end{proof}

\section{Comments}\label{sec:4}

\subsection{$\mathbf{PM}_\sigma$ versus $\mathbf{\lambda'}$, Hurewicz property and $\mathbf{AFC'}$ }\label{sec:4.1}

The relationship of classes of subsets of $\cantor$ with the Hurewicz property (and no perfect subsets, cf. Proposition \ref{examples}(1)), or of $\lambda'$-sets (cf. Proposition \ref{examples}(2)), or of $\mathbf{AFC'}$-sets (cf. Theorem \ref{AFC' sets}) to an apparently larger class of $\mathbf{PM}_\sigma$-sets seems particularly close.

 E.g., the characterizations of  $\mathbf{PM}_\sigma$-sets given in Theorem \ref{Cantor char} somewhat resemble the following characterization of sets  with the Hurewicz property, obtained by Just, Miller, Scheepers and Szeptycki \cite[Theorem 5.7]{j-m-s-s}: $A\sub \cantor$ has the Hurewicz property if and only if  for every sequence $\{K_n: n \in \NN\}$ of copies of the Cantor set in $\cantor$ disjoint from $A$ there are closed sets $F_n$ in $\cantor$ such that $A\sub \bigcup_n F_n$ and $K_m\cap F_n=\emptyset$ for each $m, n\in\NN$. In particular, if $A\sub \cantor$  has either the Hurewicz property or  $\lambda'$-property, then for every  countable set $D\sub \cantor$ disjoint from $A$ there is an  $F_{\sigma}$-set $F$ in $\cantor$ such that $A\sub F$ and $F\cap D=\emptyset$. 
 
 The following example, based on a result of Bartoszy\'nski and Shelah \cite{b-sh} and  classical ideas of Rothberger (cf. \cite{mi-1}), shows that there exists (in ZFC) a countably perfectly meager set in $\cantor$ of cardinality $\mathfrak{b}$ which lacks the latter property and thus has neither the Hurewicz nor $\lambda'$ property. It also shows that the Hurewicz and $\lambda'$ properties are not the same.

\begin{example}\label{example 1}
Inductively, one easily constructs a  subset $\{f_\alpha:\alpha< \mathfrak{b}\}$ of 
 $\baire$
 with the following properties (cf. Proposition 
 \ref{examples}(2(a))):

\begin{itemize}
	
	\item $f_\alpha$ is strictly increasing,
	\item $\alpha<\beta<\mathfrak{b}$ implies $f_\alpha<^*f_\beta$,
	\item for every $f\in\baire$ there is $\alpha<\mathfrak{b}$ with $f_\alpha\nleq^*f$.	
\end{itemize}

  By identifying each $f_\alpha$ with the characteristic function of its  range, we obtain a homeomorphic copy $A$ of $\{f_\alpha:\alpha< \mathfrak{b}\}$ in $\cantor$. 
  
  Let $B=A\cup \QQ$, where
  $\QQ$ consists of all eventually zero binary sequences.
Then:

\begin{itemize}
	\item $\{f_\alpha:\alpha< \mathfrak{b}\}$, being unbounded in $\baire$, does not have the Hurewicz property,
	
	\item  $\{f_\alpha:\alpha< \mathfrak{b}\}$ is a $\lambda'$-set in $\baire$ (cf. Proposition \ref{examples}(2)(a)),
	
	\item $B$ has the Hurewicz property and has no perfect subsets (cf. \cite[Theorem 1]{b-sh}, \cite[Theorem 2.12]{tsab} and Remark \ref{Hurewicz conjecture} below) so, by Proposition \ref{examples}(1), $B$ is countably perfectly meager in $\cantor$, 
	
	\item $A$ is countably perfectly meager in $\cantor$ as a subset of $B$,
	
	\item $A$ does not the Hurewicz property as the homeomorphic image of $\{f_\alpha:\alpha< \mathfrak{b}\}$,
	
	\item if  $F$ is any $F_{\sigma}$-set in $\cantor$ such that $A\sub F$, then $F\cap \QQ\neq\emptyset$ (since otherwise $F$ viewed as a subset of $\NN^{\NN}$ is bounded, whereas $A$ is unbounded). In particular, neither $A$ nor $B$ are $\lambda'$-sets in $\cantor$.
\end{itemize}

\end{example}

\begin{remark}\label{Hurewicz conjecture}
	It seems useful to indicate another proof of the fact that the set $B$ described in Example \ref{example 1} has the Hurewicz property. In fact, in the 
	argument below
	it is enough to assume that $\{f_\alpha:\alpha< \mathfrak{b}\}$ is any well-ordered by eventual domination and unbounded subset of $\baire$,
	$h:\NN^{\NN}\w P$ is a homeomorphism onto the set $P$ of irrationals in $I=[0,1]$,  $\QQ=I\setminus P$, $A=\{h(f_\alpha):\alpha< \mathfrak{b}\}$ and $B=A\cup \QQ$.
	
	Let us recall that the Hurewicz property of  $H\sub I$ is equivalent to the following covering property (this is the original Hurewicz's definition, cf. \cite{tsab}): 
	
	\ss 
	
	for each sequence $\U_1, \U_2,\ldots$ of open in $I$ covers of $H$,
	 there are finite subfamilies $\F_n\sub \U_n$ such that $H\sub \bigcup_n \bigcap_{m\geq n}(\bigcup \F_m)$.
\end{remark}
	
	\begin{proof}
	Let $\U_1, \U_2,\ldots$ be  open in $I$ covers of $B$
	and let $G=\bigcap_n(\bigcup \U_n)$.
	 Since $I\setminus G=P\setminus G$ is \s-compact and so is $S=h^{-1}(P\setminus G)$, there is $f\in  \NN^{\NN}$ such that $g \leq^*f$ for any $g\in S$. Let us pick $\alpha<  \mathfrak{b}$ so that 
	$f_\alpha \nleq^* f$. Then $T=\{g\in \NN^{\NN}: f_\alpha\leq^* g \}$ is an $F_\sigma$-set in $\NN^{\NN}$ disjoint from $S$ and $f_\beta \in T$, whenever $\beta\geq \alpha$.
	
	Let $L=h(T)\cup \QQ$. Then  $|B\setminus L|< \mathfrak{b}$, as $B\setminus L\sub \{h(f_\xi): \xi <\alpha\}$. Consequently, $B\setminus L$ has the Hurewicz property and since $B\setminus L\sub G$, there are finite collections
	$\F'_n\sub \U_n$ such that $B\setminus L\sub \bigcup_n \bigcap_{m\geq n}(\bigcup \F'_m)$.
	
	Let us notice that the set $L$ is $\sigma$-compact in $I$ so it also has the Hurewicz property and since $L\sub G$, we can pick finite
	$\F''_n\sub \U_n$ such that $L\sub \bigcup_n \bigcap_{m\geq n}(\bigcup \F'_m)$.
	
	Then, letting $\F_n=\F'_n\cup \F''_n$, we get finite collections $\F_n\sub \U_n$ such that $B\sub \bigcup_n \bigcap_{m\geq n}(\bigcup \F_m)$. This shows that $B$ has the Hurewicz property.
\end{proof}

\ms 

The fact that $\mathbf{AFC'}$-subsets of $\cantor$ are closely related to $\mathbf{PM}_\sigma$-sets in $\cantor$ is revealed by the following characterization of Bartoszy\'nski (private communication). We are grateful to Tomek Bartoszy\'nski for allowing us to include his result in our paper.

\begin{proposition}[T. Bartoszy\'nski]
	
For a set $A\sub \cantor$ the following are equivalent:
\begin{enumerate}
	\item $A\in \mathbf{PM}_\sigma$.
	
	\item for every  perfect subset $P$ of $\cantor$, there exists an $F_\sigma$-set $F$ in $\cantor$ such
	that $A \sub F$ and $F\cap (P+q)$ is meager in $P+q$ for every $q\in\QQ$,
		where $\QQ$ consists of all eventually zero binary sequences.
\end{enumerate}
\end{proposition}	

\begin{proof}
	
Implication $(1)\Ra (2)$ follows directly from the definition of countably perfectly meager sets.

\smallskip

To prove that $(2)\Ra (1)$, we shall appeal to condition (2) of Theorem \ref{Cantor char}.
So let   $\{K_n: n \in \NN\}$ be a sequence of Cantor sets in $\cantor$.
For each $n$ one can pick a non-empty relatively open compact set $L_n$ in $K_n$ (an intersection
of $K_n$ with a basic neighborhood in $\cantor$) and $q_n \in \QQ$ such that all points in $T_n =  L_n + q_n$  have $n$ first coordinates zero, and the sets $T_n$ are pairwise disjoint. Then the union $P$ of $\{0\}$ and the sets $T_n$ is perfect and $(K_n+q_n) \cap P$ has relatively non-empty interior in $P$ for each $n\in \NN$.

 Using (2), we fix an $F_\sigma$-set $F$ in $\cantor$ such
that $A \sub F$ and $(F+q)\cap P$ is meager in $P$ for every $q\in\QQ$.

 But now it is clear that $K_n\subseteq F$ for no $n\in\NN$, since otherwise we would have $(K_n+q_n)\cap P\sub  (F+q_n)\cap P$, contradicting the choice of $F$.

\end{proof}

The following example is a slight modification of a remarkable construction of Rec\l aw \cite{re-2} and shows that, at least consistently, there exists  a countably perfectly meager in $\cantor$, in fact a $\lambda'$-set in $\cantor$, which is not  perfectly meager in the transitive sense (a similar construction was also used by Weiss \cite[Theorem 3]{w} in his proof that the existence of a universally meager set in $\cantor$ of cardinality of the continuum implies that there is also one which is not perfectly meager in the transitive sense). This is also yet another (cf. Theorem \ref{lambda}) strengthening of the result of Nowik, Scheepers and Weiss \cite{n-s-w} that under the Continuum Hypothesis there is a $\lambda$-set in $\cantor$ which is not  perfectly meager in the transitive sense.

 \begin{example}\label{example 2}
Let us assume that there exists a $\lambda'$-set in $\cantor$ of cardinality of the continuum.

Let $C,\ D$ be disjoint copies of the Cantor  set in $\cantor$ such that 
\begin{enumerate}
	\item[(1)] the operation $+$ of addition is a homeomorphism between $C\times D$  and $C+D$ (cf. \cite{re-2}).
\end{enumerate}

Let $T\sub C$ be a $\lambda'$-set in $C$ of cardinality of the continuum.

Arguing as in part (A) of the proof of Theorem \ref{main}, we obtain a set $H\sub T\times D$ intersecting each vertical section $\{t\}\times D$, $t\in T$, in a singleton, such that
\begin{enumerate}
	\item[(2)]  each $F_\sigma$-set in $C\times D$ containing $H$ contains also $\{t\}\times D$ for some $t\in T$.
\end{enumerate}

Now, $T$ being a  $\lambda'$-set in $C$, one readily checks that $H$ is a $\lambda'$-set in $C\times D$ (cf. the proof of Theorem \ref{lambda}(3)). It follows that, cf. (1),  if we let $Y=+(H)$, then $Y$ is a $\lambda'$-set in $C+ D$ and hence also in $\cantor$.

On the other hand, $Y$ is not perfectly meager in the transitive sense. Indeed,  if $F$ is an arbitrary $F_\sigma$-set in $\cantor$ with $Y\sub F$ and we let $E= (+)^{-1}(F\cap (C+ D))$, then $E$ is an $F_\sigma$-set in $C\times D$ containing $H$ so it also contains, cf. (2), $\{t\}\times D$ for some $t\in T$. Consequently, $t+D\sub F$, completing the proof.

\end{example}

\subsection{Remarks on related \s-ideals}\label{sec:4.2}
In this section  $\F$ always denotes a countable (possibly finite but  non-empty) collection of perfect sets in a perfect Polish space $X$.

If $I$ is a \s-ideal of subsets of $X$, then by $I^*$ we  denote the \s-ideal generated by the
closed subsets of $X$ which belong to $I$. 
 
\subsubsection{The \s-ideals $MGR(\F)$}\label{MGR(F)} 
 Following Kechris and Solecki \cite{ke-so} let us put
$$
MGR(\F)=\{A\sub X: A\cap P \hbox{ is meager in $P$ for every } P\in\F  \}.
$$

Let us note that $MGR(\F)$ is a \s-ideal on $X$ 
 generated by Borel, in fact $F_{\sigma\delta}$-sets, and fulfills the c.c.c. Moreover, the quotient Boolean algebra $Bor(X)/(MGR(\F)\cap Bor(X))$ is isomorphic to the Cohen algebra, cf. \cite{b-2}. Clearly, the intersection of all \s-ideals of the form $MGR(\F)$ is precisely the \s-ideal $\mathbf{PM}$.

 On the other hand, by \cite[Theorem 2.1]{z-1}, the \s-ideal $\mathbf{UM}$ is the intersection of all \s-ideals $I$ on $X$ such that the quotient Boolean algebra $Bor(X)/(I\cap Bor(X))$ is isomorphic to the Cohen algebra. Any such $I$  is precisely of the form $\M(X,\tau)$, by which we denote the \s-ideal consisting of meager sets 
with respect to a perfect  Polish topology $\tau$ 
on $X$  giving the original Borel structure of $X$.

  \subsubsection{The \s-ideals $MGR^*(\F)$}\label{MGR*(F)} 
  
    By the definition (see the beginning of Section \ref{sec:4.2}), the \s-ideal $MGR^*(\F)$ consists of such sets $A\sub X$ that there exists an $F_\sigma$-set $F$ in $X$ with $A\sub F$ and such that $F\cap P$ is meager in $P$ for every  $P\in\F$. It is the \s-ideal $MGR^*(\F)$ (not the \s-ideal $MGR(\F)$ erroneously  employed in \cite[Theorem 7, $(1)\Ra (3)$]{b-2})   
  that is relevant to the definition of countably perfectly meager sets. Indeed, the intersection of all \s-ideals of the form $MGR^*(\F)$ is precisely the \s-ideal $\mathbf{PM}_\sigma$. 
  
  It turns out that not all of the \s-ideals of the form 
  $MGR^*(\F)$ fulfill the c.c.c. Let us elaborate on this  a little further with the help of a theory developed by Kechris and Solecki \cite{ke-so}.
  
  \begin{proposition}\label{char 4.1}\hfill\null
  	\begin{enumerate}
  		
  		\item If  every non-empty open set $U$ in $X$ contains a nowhere dense set $P\in\F$, then the \s-ideal $MGR^*(\F)$ does not fulfill the c.c.c.
  		
  		\item The intersection of all \s-ideals in $X$ generated by closed sets which fulfill the c.c.c. is precisely the \s-ideal $\mathbf{PM}$.
  		
  		\item The intersection of all  \s-ideals of the form $MGR^*(\F)$ which fulfill the c.c.c. is precisely the \s-ideal $\mathbf{PM}$.
  		
  		\item The intersection of all \s-ideals of the form $MGR^*(\F)$ which do not fulfill the c.c.c. is precisely the \s-ideal $\mathbf{PM}_\sigma$.
  		  		
  	\end{enumerate}
  \end{proposition} 

\begin{proof}
	
	(1). This immediately follows from \cite[Lemma 9]{ke-so}.
	
	(2) and (3). First, let $A\in \mathbf{PM}$ and let $I$ be a \s-ideal in $X$ which is generated by closed sets and fulfills the c.c.c. Then, by \cite[Theorem 3]{ke-so}, $I$ is of the form $MGR(\F)$ for a countable family of perfect subsets of $X$. Consequently, $A\in I$.
	
	For the other direction, assume that $A\in MGR^*(\F)$ for every $\F$ such that the \s-ideal $MGR^*(\F)$ fulfills the c.c.c. Suppose that $P$ is an arbitrary perfect subset of $X$ and let $\F=\{P\}$. Then we have $MGR^*(\F)=MGR(\F)$, so  the \s-ideal $MGR^*(\F)$ fulfills the c.c.c. Consequently, $A\in MGR^*(\F)$ which just means that $A\cap P$ is meager in $P$, completing the proof that $A\in \mathbf{PM}$. 
	
	(4). Let us assume that $A\in MGR^*(\F)$ for every $\F$ such that the \s-ideal $MGR^*(\F)$ does not fulfill the c.c.c.  To prove that $A\in \mathbf{PM}_\sigma$, let $\F$ be an arbitrary countable family of perfect subsets of $X$. By extending $\F$, if necessary, we may assume that every non-empty open set $U$ in $X$ contains a nowhere dense set $P\in\F$. By part (1), the \s-ideal $MGR^*(\F)$ does not fulfill the c.c.c. Consequently, $A\in MGR^*(\F)$ and we are done.

\end{proof}
	
	\subsubsection{The \s-ideals $\M^*(X,\tau)$}\label{M^*} 
In this subsection  $\tau$ always denotes  a perfect  Polish topology
on $X$  giving the original Borel structure of $X$. Recall that $\M(X,\tau)$ is the \s-ideal  of meager sets 
with respect to  $\tau$ and $\M^*(X,\tau)$ consists of such $A\sub X$ that there exists an $F_\sigma$-set $F$ in $X$ (with the original Polish topology) with $A\sub F$ and  $F\in\M(X,\tau)$. 

It turns out that we can characterize perfectly meager and  countably perfectly meager sets with the help of the \s-ideals $\M^*(X,\tau)$ in an analogous way to their characterizations in terms of the \s-ideals $MGR^*(\F)$ (cf. \ref{char 4.1}).  

 \begin{proposition}\label{char 4.2}\hfill\null
	\begin{enumerate}

		\item The intersection of all \s-ideals of the form $\M^*(X,\tau)$ which fulfill the c.c.c. is precisely the \s-ideal $\mathbf{PM}$.
		
		\item The intersection of all \s-ideals of the form $\M^*(X,\tau)$ is precisely the \s-ideal $\mathbf{PM}_\sigma$. Moreover, for a set $A$ to be countably perfectly meager it is enough to belong to all \s-ideals of the form $\M^*(X,\tau)$ which do not fulfill the c.c.c.
		
	\end{enumerate}
\end{proposition} 

\begin{proof}
	
	(1). This follows directly from points (2) and (3) of Proposition \ref{char 4.1}, since every \s-ideal of the form $MGR^*(\F)$ is also of the form $\M^*(X,\tau)$ which in turn is generated by closed sets. 
	
	(2). Let us assume that  $A\in \mathbf{PM}_\sigma$ and let $\{U_n:n\in\NN\}$ be a basis of a perfect  Polish topology $\tau$ on $X$ giving the original Borel structure of $X$. We assume with no loss of generality that $\tau$ extends the original Polish topology on $X$.
	  For each $n$ let us pick a Cantor set $K_n\sub U_n$. We complete the argument exactly as in the proof of implication $(1)\Ra (2)$ in Theorem \ref{Baire sp char}.
		 
	 For the other direction we again use the fact that every \s-ideal of the form $MGR^*(\F)$ is also of the form $\M^*(X,\tau)$ and apply point (4) of Proposition \ref{char 4.1}.

\end{proof}

\bibliographystyle{amsplain}

\end{document}